\documentclass[12pt,thmsa]{article}
\usepackage{amsmath}
\usepackage{amssymb}

\usepackage[dvips]{graphicx}
\usepackage{amssymb}
\usepackage{amsmath}
\date{}

\newtheorem{theorem}{Theorem}[section]

\newtheorem{lemma}[theorem]{Lemma}

\def\1{\bf 1}

\begin{document}
\date{}
\author{ Volodymyr Koshmanenko$^1$ and Igor Samoilenko$^2$}
\title{\bf The  conflict triad dynamical system} \maketitle

\begin{abstract}
A  dynamical   model of the natural conflict triad  is investigated.
The conflict interacting  substances of the  triad are: some
biological population, a living resource, and a negative factor
(e.g., infection diseases). We suppose  that each substance is
multi-component. The main  coexistence phases for substances are
established: the equilibrium  point (stable state), the local cyclic
orbits (attractors), the global periodic oscillating trajectories,
and the evolution close to chaotic. The bifurcation points and
obvious thresholds between phases  are exhibited in the computer
simulations.
\end{abstract}

$^1${\ Institute of Mathematics,  Tereshchenkivs'ka str. 3, Kyiv
01601 Ukraine} \ {E-mail: kosh@imath.kiev.ua}

$^2${\ Institute of Mathematics,  Tereshchenkivs'ka str. 3, Kyiv
01601 Ukraine} \ {E-mail: isamoil@imath.kiev.ua}



\section{Introduction} We use the term   {\it
conflict triad} for notation of a physical system consisting of
three conflict substances (denoted by ${\bf P, R, Q}$) which exist
in the same space $\Omega$.  Each
 substance at  initial moment of
time is presented by a real value $P, R$, or $Q$ in accordance
with its amount characteristic. The whole system is complex since
   every substance contains a family of components
   distributed along  regions $\Omega_i$ which compose the existence
   space: $\Omega=\bigcup_{i=1}^{n}\Omega_i, \ 2\leq n< \infty.$
So one can think that  $ {\bf P, R, Q}$  are  vectors with
non-negative coordinates:
 $$  {\bf P}=(P_1,...,P_n), \  {\bf R}=(R_1,...,R_n), \  {\bf Q}=(Q_1,...,Q_n),
  \ P_i, R_i, Q_i\geq 0, \ i=1,...,n,$$
$$P= P_1+ \cdots + P_n ,\ R= R_1+ \cdots + R_n, \ Q= Q_1+ \cdots +
Q_n.$$

We study the evolution of the system at discrete time: $$\{{\bf
P}^{N}, {\bf R}^{N}, {\bf Q}^{N} \}
 \stackrel{\divideontimes
}{\longrightarrow } \{{\bf  P}^{N+1}, {\bf R}^{N+1}, {\bf Q}^{N+1}
\}, \ \ N=0,1...,$$
 where ${\bf P}^{0}={\bf P}, \ {\bf R}^{0}= {\bf R}, {\bf Q}^{0}={\bf Q}
 $ and the map    $ \stackrel{\divideontimes
}{\longrightarrow } $ is defined by the law of conflict interaction.
In general this law is unknown. Our definition of the transformation
$\divideontimes$  (see below formulae (\ref{gl}), (\ref{reg}))
starts with the famous Lotka-Volterra predator-prey approach (see
e.g. \cite{HoSi1,HoSi2,Mu1,MM}) and is based on heuristic
understanding of the physical essence of the substances ${\bf P, R,
Q}$.

The schematic picture of mutual dependencies between ${\bf P, R,
Q}$ may be  presented by the following diagram:

$\begin{picture}(10,6)(1,15) \put (8,19){\circle{10}} \put
(4,15){\circle{10}} \put (12,15){\circle{10}} \put(7.9,18.8){P}
\put(3.9,14.8){Q} \put(11.9,14.8){R} \put (4,16){\vector(1,1){3}}
\put (7.8,18.2){\vector(-1,-1){3}} \put
(11.2,15.2){\vector(-1,1){2.8}} \put (9,18.8){\vector(1,-1){2.8}}
\put (6,15.3){\vector(1,0){4}} \put (10,14.2){\vector(-1,0){4}}
\put(5,17.5){--} \put(6.5,16.5){+} \put(9.2,16.5){+}
\put(10.5,17.5){--} \put(7.9,15.5){--} \put(7.9,13.7){--}
\end{picture}$

\ \ \ \ \

\ \ \ \ \ \

\ \ \ \ \ \ \ \ \

where each signed arrow  means  a direction of positive or negative
influence   one of substances at another. The positive dependence
involves a certain growth and negative -- decreasing of quantitative
characteristic of the corresponding   specie.

At this diagram  we have three types of interaction connections: \
two variants of the  \ "plus--minus" \  interactions  and a single
version of the \ "minus--minus" \ interdependence.

The
 "plus--minus" \ type of interaction is in fact the modified stochastic  version
 of the  \  predator -- prey
 model (c.f. with \cite{E.R}).
 Here we use two different variants of
this type. First, the conflict  interdependence between behavior of
some biological population playing a role of  predator and the
dynamical changes of the vital resource as a prey. Second, the
nonlinear relations between the local density of biological species
(as some kind of a prey) and the expansion processes of viral
infections, like a certain predator.

The "minus--minus" \ interaction is essentially different. This type
describes  the conflict competition of mutually alternative
substances like  virus infection diseases and medicine. In the paper
we consider this couple slightly  wider, we treat the vital
resources as the opposite substance to the viral infection and admit
its quantitative decreasing under destroying actions of the latter.

 In more details, one can understand a sequence $P^N_i, N=1,2,...$
as the quantitative evolution of some  biological specie in a fixed
region $\Omega_i$ under  influence of two opposite factors. Namely,
the population dynamics of $P^N_i$ is determined by the positive
dependence produced by the vital  resource environment from the side
${\bf R}$ that
 provides some growth  of $P^N_i$, and  by the negative influence
generated by  ${\bf Q}$. The possible viral infection caused by the
substance ${\bf Q}$ leads to the quantitative losses for $P^N_i$.
Thus, at the moment
 $N+1$ the value
  $P^{N+1}_i$ has to be proportional to
amount of the biological spices  in region $\Omega_i$ at the
previous moment, it should grow at increasing of $R^N_i$, but to
fall at increasing of $Q^N_i$.

In turn, under assumption the vital resource is not exhaustible
globally, the regional changes of its amount $R^N_i$ are dependent
naturally from the local (=regional) intensity of their using as a
source for existence of biological species, and in addition from the
negative influences of  ${\bf Q}$ which has the alternative nature.
Finally, the evolution of coordinates $Q^N_i$ describes the
dynamical picture for behavior of  threat concentrations for
existence of biological species inside every separate region
$\Omega_i$. So the local expansion of some infection \ $Q^{N+1}_i$
grows  at increasing of $P^N_i$ and decreasing of $R^N_i$.

The explicit  formulae of all mutual interdependence will be defined
in  Section 3.

One can associate the triple of substances   ${\bf P, R, Q}$  with
the philosophical triad: mankind, good, and evil or, in accordance
with mythology, with  flora and fauna of Earth surrounded by water
and fire as a positive and negative resources. More specifically,
$P^N$ may mark the quantitative global population of some biological
specie (humanity),  $ R^N$ -- the variation  of vital  resources in
the current environment, and $ Q^N$ corresponds to the dynamics of
threatening  concentrations associated with various negative
phenomena that may influence both at the biological population $P^N$
and the resource environment $ R^N$ too. In particular, coordinates
$P^N_i, R^N_i, Q^N_i$ describe the  local evolution of components of
triple indicated substances in the fixed region $\Omega_i$.

The general scheme of  the dynamical system studied here is complex
enough and it is hardly to expect for establishment of abstract
mathematical results (theorems) concerning the global behavior.
Nevertheless, we got some exact results which are true, in
particular, for separate two-sided  links inside of the conflict
triad, i.e., for couple interaction between fixed substances (see
Section 2). We study the whole complex system by using the series of
computer models in Section 3. Their analysis exhibits remarkable
features which would be useful in applications.

Two intrinsic questions arise. Whether the construction of conflict
triad is well-defined? If so, whether the  dynamical picture of the
substances ${\bf P, R, Q}$ coincides with our intuitive pattern  of
the conflict triad coming from natural practical experience?

The first important result of this work is the setting of
mathematical consistency  of our constructions. All the
interrelations  used in  the dynamical system of the conflict triad
do not lead to collapse. More exactly, it is shown that under a
certain choice  of additional parameters of the model (this choice
is the non-trivial  problem itself)  the whole complex system does
not decay (destroy themselves). This means that the values of all
coordinates are changed in the physically reasonable  scopes.

The second  important result of our work  is that various computer
models that realize the conflict triad expose the typical features
for the complex systems behavior. Namely, analysis of the values
$P_i^{N}, R_i^{N}, Q_i^{N}$ in concrete models demonstrate the
distinctive properties that are observed in natural conditions where
more than two substances take part in confrontations. Among such
properties, in particular, we observe  the  presence of distinct
phases for existence of the system: \ the phase  of dynamic
equilibrium (the equilibrium state), the phase of periodic
oscillations which contains the wide spectrum of trajectories that
have tendency to approximate the cyclic attractors, the presence of
areas for bifurcation points and thresholds between different
phases, finally we meet the class of evolution close to  chaotic
(which is called the quasi-chaotic).

 We note that  mathematical basis  of our constructions
  is a concept of the conflict dynamical system that has been developed
 in  \cite{KoTC, KoTC1, KoDo, BKoDo, AKS, ABodK, KoKh}, (c.f. with
 \cite{SaTa}).
 Among other deep theoretical researches that influenced on our work
  we refer \cite{ABG,BGGLP},\cite{Ku} - \cite{E.R}.

The most close  to our conflict triad is the well-known SIR-model
and its  different variants  that describe the dynamic  of epidemic
infections (for  details see, for example, \cite{Mu1,Mu2}). In the
SIR-model the idea of conflict triad is actually presented but in an
implicit  form. So, all the population of some biological specie is
divided into three groups: $S$ is the amount of persons favorable to
the infection (susceptible), $I$ -- infected and able to carry the
infection, and $R$ -- which are under cover or recovered already.
The complicated SIR-model contains the additional group $E$ of
those, who has an infection in the hidden (latent) form. The
evolution of the above mentioned groups in time is determined by
relatively simple but nonlinear equations:
$$\dot{S}=-rSI, \dot{I}=rSI-aI, \dot{R}=aI,$$ where coefficients
$r, a>0.$ It is assumed that $S(t)+I(t)+R(t)=const$ and since that
$\dot{S}+\dot{I}+\dot{R}=0.$ The main problem is to research the
dynamic of distribution for the infection $I(t)$, in particular,
when $I(t)$ increases or decreases depending on coefficients $r, a$
and starting values $S(0), I(0), R(0).$ \  There exist a number of
publications (see references in \cite{Mu1,Mu2,Ep})  where the
SIR-model is improved or modified in accordance with concrete
specifics of biological species, view of epidemic infections, and
conditions of its expansion in a certain environment.

In comparison with the  SIR-model, our construction of the
conflict triad is substantially more perfect in two principal
aspects.

First, we use a partition of whole  existence space into
 family of finitely many regions. Of course, it corresponds to
 the widely  observed environment phenomenon of nature:  the existence space
 always is separated into bounded domains. That is why
the computer models of the system are able, in particular, to
describe the expansion of disease infections of biological
population separately in every region of existence and to define,
for example, the most safe regions.

Another new important step in our work is that we put as a basic the
probabilistic (statistical) law of interaction between conflict
substances. We consider that statistical interdependence reflect
more deep connections between opponents any nature. Thus we may
describe not only quantitative changes of populations but also its
redistribution along  regions of existence caused by the conflict
interaction.

Shortly, in contrast to the above mentioned  SIR-model our method of
the construction of the conflict triad uses a partition of the
existence space into regions  and application
 of  statistical
formulas for description of the dynamical picture.

Finally, it is especially important, we assume that each substance
(opponent) is {\it a priori}  non-annihilating. This  principle is
putting directly in  formulas (\ref{gl}), (\ref{reg}) that govern
the conflict interactions. So,  any of  substances ${\bf P, R, Q}$
can not annihilate another one, but has to reach the compromise
 equilibrium  state or  migrate along regions by some law.

\section{Dynamics of the bilateral couple interactions}

According to equations (\ref{gl}), (\ref{reg}) (see below) the
interactions between substances $\bf P, R, Q$ in  the conflict triad
are nonlinear and rather complex. But the separate bilateral couple
interdependency  admits rigorous  enough mathematical analysis. In
this section we state two results for situations marked as "plus --
minus" \ and "minus -- minus" \ models.

In the first bilateral  "plus -- minus" \ model we analyze the
dynamics of changes for substances  $\bf P$ and  $\bf R$. Plus means
that $\bf P$ is positively influenced by  $\bf R$, and minus -- that
$\bf R$ is negatively influenced by $\bf P$. That is a typical
predator -- prey situation (see for example \cite{MM}). However we
consider different, more specific  functional dependence between
opponents.

 Let $P^N$, \ $N=0,1,...$, denote the amount
 evolution  of  $\bf P$  (some fixed
biological specie)   at discrete time. As soon as the space $\Omega$
is divided onto family of separate regions \ \
$\Omega=\bigcup_{i=1}^{n}\Omega_i, \ 2\leq n<\infty,$ \ \  the
complete amount $P^N$ at every time moment is certainly distributed
along these regions: $ \ P^N= P^N_1+\cdots +P^N_n$. For
simplification of the problem, we  assume that the complete amount
of substance at the whole territory $\Omega$ is permanent \ $P^N= P=
 const.$ \ It means that the average rate of growth and the decay
rate of the substance $\bf P$ population (the rates of birth and
deaths) are independent by the time, i.e. are in the dynamical
balance. In fact this occurs  often enough, at least locally at
certain periods of development of a complex physical system. Thus,
in the such simplified model we analyze only the redistribution
dynamic for $ P_i^N$ along the regions $\Omega_i$. In the model of
conflict triad the redistribution of values $ P_i^N$ is stipulated
by two factors: the positive influence of the substance $\bf R$ \
(the population of  $ P_i^N$ grows using the vital resource
environment), \ and the negative influence from the side of $\bf Q$
(for example, threats,  infection diseases  cause some  decay). In a
situation, where the global influence of the last factor is
insignificant and negligible, the bilateral coupled interactions
between regional values of $P^N$, $R^N$ are essentially simplified
in the mathematical sense. In turn, the dependence of $R$ from $P$
is purely negative (the amount of the substance $\bf R$ is "burned"
\ as a vital resource for ${\bf P}$). However we assume  the global
amount of $R$ in $\Omega$ is stable $R^N= R= const$. This stability
condition means that $\bf R$ is continuously supplied due to the
external source (such as the Sun) that now is not examined.
Nevertheless at each time moment   the regional distribution of the
vital resource
 $ R^N= R^N_1+\cdots +R^N_n$  is changed according to the
appropriate law.

Thus the simplified problem is to study the redistribution dynamic
 for values $ P_i^N$, $R^N_i$ along regions
$\Omega_i$ produced by a certain  "plus -- minus" \ interaction
between substances ${\bf P}$ and $\bf R$ under assumption of the
global  amount stability for $P$, $R$ in  $\Omega$ and their
independence with the third substance.

The proper   conflict dynamical system has the view:
\begin{equation}\label{dscPR} \{P_i^{N}, R_i^{N}\}
\stackrel{\divideontimes }{\longrightarrow }\{P_i^{N+1},
R_i^{N+1}\}, \ \ P_i^{0}=P_i, R_i^{0}= R_i.
\end{equation}  Surely it is considerably simpler as compared to the behavior
of all complex system. After the explicit definition of the conflict
transformation
 $\divideontimes$ (see below (\ref{pm})) we are able to fulfill
 the detailed enough analysis of system (\ref{dscPR}). The main results are
 stated in Theorem 1 which gives complete enough description of behavior of this system.

\ \ \ \ \

Let us introduce the simplest variant of concrete formulas of the
\ "plus -- minus" \ interaction between substances ${\bf P}$ and
$\bf R$ . We write down these formulas in terms of coordinates of
stochastic vectors with a unite  $l_1$-norm:
\begin{equation}\label{PR} {\bf p}^N=(p^N_1...p^N_i,...,p^N_n),
\ \
 \ \ {\bf r}^N=(r^N_1...,r^N_i,...,r^N_n), \ \ \
p^N_i:=\frac{P_i^N}{P}, \ \ r^N_i:=\frac{ R_i^N}{R},\end{equation}
where, we remind, $P, R$ are the  amount characteristics of
substances ${\bf P}$, $\bf R$ in the whole space $\Omega$.
 Namely, at $N+1$ step the coordinates
$p^{N+1}_i$ and $r^{N+1}_i$ are determined iteratively by the rule
\begin{equation}\label{pm}
p_i^{N+1}=\frac{p_i^{N}(1+r_i^{N})}{z_{p,r}^{N}}, \ \
r_i^{N+1}=\frac{r_i^{N}(1-p_i^{N})}{z_{r,p}^{N}}, \ \ N=0,1...
\end{equation} where the normative denominators $$z_{p,r}^{N}=1+({\bf
p}^{N},{\bf r}^{N}), \ \ z_{r,p}^{N}=1-({\bf p}^{N},{\bf r}^{N})$$
($(\cdot,\cdot)$ denotes the inner product in ${\mathbb R}^n$) \
ensure that vectors ${\bf p}^{N+1}, {\bf r}^{N+1}$ remain
stochastic.

We note that different signs in the numerators of formulas
(\ref{pm}) just determine the essence  of the \ "plus -- minus" \
model. Plus means that $P_i^{N+1}$ increases depending on a value
$R_i^N$, and minus supplies the  decay of $R_i^{N+1}$ at increasing
$P_i^{N}$. Of course, one have to make these interpretations after
re-normalizing inverse to (\ref{PR}). In fact, the model is
well-defined by virtue of
 the normative denominators.

\ \ \ \ \ \ \
\
\
\
\
\

\begin{theorem} \label{Th1} ("Plus-minus" \ model) \ The conflict
dynamical system
 (\ref{dscPR}) given  by the formulae  (\ref{pm}) has three
typical phases of behavior.

The first phase (equilibrium point) is determined by the uniform
distribution:
  \begin{equation}
 \label{Iph} \forall N \ \ P^N_i=P/n, \ \ R^N_i=R/n;
  \ \ {\bf p}=(1/n...,1/n), \ \ {\bf r}=(1/n...,1/n).
\end{equation}
However this equilibrium state  is unstable.

The second phase (existence of the limiting fixed points) is
stipulated  by the condition: at least for one $i$,
 \begin{equation}
 \label{IIph}p_i=0  \ \
r_i\neq 0.
\end{equation}
 In this case the  system  trajectory  converges to the limiting stable
 state,
   \begin{equation}\label{pr}
{\bf p}^{\infty} = \lim_{N\to \infty}{\bf p}^{N}, \ \ {\bf
r}^{\infty} =  \lim_{N\to \infty}{\bf r}^{N}, \ \ {\bf p}^\infty
\perp {\bf r}^\infty,
\end{equation}
which is invariant with respect to the conflict interaction.

 The
third typical phase (the quasi-chaotic behavior)
  occurs under the starting conditions
 \begin{equation}  {\bf p}\neq {\bf r}, \ \ \forall i, p_i \neq
0 \  r_i \neq 0.\end{equation} At this phase each of coordinates
$p_i^N, r_i^N \ i=1...,n$ oscillates between zero and one, in
general case without any regular law. \end{theorem}
 \ \ \ \ \ \

{\it Proof.} \ In the case of uniform starting distributions of
${\bf P}$, $\bf R$ \ along regions, $P_i=P/n, \ R_i=R/n, \forall i,$
all coordinates of the stochastic vectors ${\bf p}, {\bf r}$ are
equal: $p_i=r_i=1/n$. Then it is easy to find that $$({\bf p},{\bf
r})=1/n=z_{p,r}=z_{r,p}.$$ Therefore due to (\ref{pm}) for all $N$
we have $p_i^N=r_i^N=1/n$. It proves that   uniform starting
distributions define the equilibrium point for dynamical system
(\ref{dscPR}). This state is unstable. An arbitrary small deviation
$\varepsilon>0$ of any coordinate $p_i$ or $r_i$ from $1/n$  leads
in time to greater deviations (see Lemma 1 below).

Let us prove  (\ref{pr}) under condition
  (\ref{IIph}). Here we introduce the value $$\theta^N:=\sum_i p_i^Nr_i^N=({\bf
p}^N,{\bf r}^N)$$ and call it the {\it conflict index} for dynamical
system at moment $N$.

Without loss of generality we assume that
 (\ref{IIph}) is fulfilled for only single coordinate and all other
ones,  $r_i, p_{k\neq i}, r_{k\neq i}$ are nonzero. Then, it
follows from
 (\ref{pm}), we have to show
 that there exist the limiting vectors ${\bf
p}^{\infty},  \ {\bf r}^{\infty}$ with coordinates:
$$p_i^\infty=0, r_i^\infty=1, \
  \ p^\infty_k\geq 0, r^\infty_k= 0, \ \ k\neq i.$$
Indeed, if   $p_i=0$ and
 $r_i\neq 0$, then due to $0< \theta\equiv \theta^0 < 1$ the sequence
   $r_i^N=r_i^{N-1}/(1-\theta^{N-1}), \ N=1,2,...$ monotonically increases.
Thus, since $r_i^N<1$, there exists the limit $r_i^\infty=\lim_{N\to
\infty} r_i^N$. That is, because
$r_i^N=r_i\cdot\prod_{l=1}^N1/(1-\theta^l)< 1$, the convergence of
$r_i^N $ implies with necessity that the  conflict index $\theta^N$
monotonically decreases. In fact $\theta^N\to 0$. To see
$r_i^\infty=1 $ we prove that all $r_{k\neq i}^N \rightarrow 0, \
N\to \infty$. Assume the opposite, i.e., that there exists at least
single coordinate $r_{k}^N, \ k\neq i $ which does not converge to
zero. Then, due to  the  conflict index  $\theta^N\to 0$, the
coordinate $p_{k}^N$ have to come to zero. By  (\ref{pm}) we have
$$p_{k}^{N+1}=p_{k}^N\cdot (1+r_{k}^N)/(1+\theta^N) \rightarrow 0,$$
that is possible if only $r_{k}^N < \theta^N$, thus it is the
contradiction.  So $r_{k}^N\to 0$. And  by the same reason all other
coordinates $r_{k\neq i}^N$ converge to zero too. Therefore
$r_i^\infty= 1$. By similar way we obtain the existence of the
limiting coordinates $p^\infty_k, k\neq i$, which may take any
non-zero values given unit in a sum. It is evident due to
$\theta^\infty = 0$ that $ {\bf p}^\infty \perp {\bf r}^\infty$ and
therefore in the limiting state the system reaches the stable
equilibrium.

Let us consider the case of the third phase. We assume that vectors
$ {\bf p}, {\bf r}$ are different and all starting coordinates are
non-zero. Then $0< \theta <1$,  the  vectors
 ${\bf p}, {\bf r}$ are not orthogonal. In particular,
 we exclude that $p_i=r_i=1$ for some $i$.
Let us show that in such a case all coordinates oscillate inside
open interval between zero and unit (without any evident pattern).

Our argumentation is based entirely on formulae (\ref{pm}). Take
any couple of coordinates $0\neq p_i\neq r_i \neq 0$. Assume that
in the starting moment the inequalities
\begin{equation} \label{e1}
r_i<p_i<\theta \end{equation} hold. We shall show that all other
possible inequalities between  $r_i^N, p_i^N, \theta^N$ will
successively appear  at some moments of time. In the next we use the
evident fact about the qualitative behavior of values
 $\theta^N$, \ $p_i^N, r_i^N$. Namely, by the definition of
the  conflict index,
   $\theta^N=\sum_i p_i^N \cdot r_i^N,$ its small changes are slower
   than changes of
   any fixed coordinate $ p_i^N, r_i^N$. Indeed, a finite sum of the differential
    products    $\Delta p_i \Delta r_i$ has the second power of smallness with respect
    to  any separate differential $\Delta p_i$ or $\Delta r_i$.
     In what
follows we omit the time subscript $N$.

Directly from (\ref{e1}) by (\ref{pm}) it follows that $p_i$
decreases and $r_i$ increases. This leads with necessity to the
transformation (\ref{e1}) into the inequality
\begin{equation} \label{e2}
p_i<r_i<\theta. \end{equation} Using again formulae (\ref{pm}) we
check that the coordinate $p_i$ will be else decreasing, and $r_i$
-- increasing. This is continued till the moment when instead
(\ref{e2}) appears the inequality
\begin{equation} \label{e3}
p_i<\theta<r_i. \end{equation} In turn, again by formulae (\ref{pm})
the coordinate $p_i$ begins to grow due to (\ref{e3}), although
$r_i$ is still increasing. This tendency continue till the moment
when the inequalities  (\ref{e3}) change at
\begin{equation} \label{e4}
\theta<p_i<r_i. \end{equation}

We note,  inequalities (\ref{e4}) do not imply  the convergent to
zero of  the conflict index, i.e., the vectors ${\bf p}, {\bf r}$
could not become orthogonal.

By   (\ref{e3}) the coordinate
 $p_i$ grows. It continues to grow after coming to (\ref{e4}), but
  $r_i$ begins to decay. This leads to the new inequalities
\begin{equation} \label{e5}
\theta<r_i<p_i. \end{equation} In turn, it follows that $p_i$
still grows, and $r_i$ continues to decay. On this way the
inequalities
\begin{equation} \label{e6}
r_i<\theta<p_i \end{equation} appear. The latter produce the
 decaying  of  $p_i$ and $r_i$ also decreases as at previous period.
 But these changes continue only till the moment when
  (\ref{e6}) is replaced by the starting inequalities (\ref{e1}).

 We note that the sequence of
transformations from (\ref{e1}) to (\ref{e6}) are ordered and
sometimes equalities may appear but by (\ref{pm})  they pass into
inequalities immediately at the next moment. Thus the full cycle
of all possible  inequalities between $\theta, r_i, p_i$ was
realized.

By the way we observe that no one coordinate  $r_i, p_i$, as well
as the conflict index $\theta$ could not go to zero closely. In
particular, for example, if  (\ref{e1}) takes place then  $r_i$
begins to grow. In general, as soon as some coordinate becomes
smallest it begins to increase  with necessity that follows from
(\ref{pm}).

It is not hard to see that for different couples $r_i, p_i$ \ each
inequality from (\ref{e1}) -- (\ref{e6}) is fulfilled at various
moments of time and the corresponding transformations are not
synchronous. That is why the existence of regular cyclic
oscillations for the all complex system in the considered situation
is questionable. Therefore in this phase orbits of the system are
similar to some kind of the quasi-chaotic behavior.

If one or several coordinates    $r_{i_1}=\cdots =r_{i_m}=0,
 1\leq m<n$,  but all $p_{i} \neq 0, \ i=1,...,n$, then it is easy to see
 that  $p_{i_1}^N,...,p_{i_m}^N  \rightarrow 0$, and residual
 coordinates have the  quasi-chaotic behavior.

The proof will be completed if we prove that the equilibrium state
(\ref{Iph}) is not a one-point attractor (see Lemma \ref{L} bellow).

\ \ \ \ $\Box$

\ \ \ \ \ \ \

\begin{lemma} \label{L} An arbitrary small deviation from  the equilibrium state  (\ref{Iph})
of the conflict dynamical system
 (\ref{dscPR}): $$p_i=1/n \to
p_{і,\varepsilon}=1/n+\varepsilon_i, \ r_i=1/n \to
r_{і,\delta}=1/n+\delta_i, \ \ \sum_i\varepsilon_i=0, \ \sum_i
\delta_i=0$$ automatically leads to the greater deviation at least
for some coordinates:$$\varepsilon_i \to
\varepsilon_i'>\varepsilon_i \ \delta_i \to \delta_i'>\delta_i.$$
\end{lemma}

{\it Proof} \ is not trivial and here we present only its sketch.
We need to analyze the dependence  of values $p_i^N=r_i^N$ from
fixed deviations as $N\to \infty$. The linearization of
 (\ref{pm}) shows that the terms of the first order by $\varepsilon_i,
 \delta_i$ expose the following changes of deviations
$$\varepsilon_i \to \varepsilon_i'=\varepsilon_i+1/(n+1)\delta_i,
\ \delta_i \to \delta_i'=\delta_i-1/(n-1)\varepsilon_i.$$ These
formulae are defined by the strictly positive definite matrix $t$
with elements $t_{11}=t_{22}=1, \ t_{12}=1/(n-1), \ t_{21}=
-1/(n-1)$ \ which does not depend on the starting deviation. By
this reason the iteration of formulas
 (\ref{pm}) produce new  deviations increasing with $N\to \infty$,
  that one may check
 directly. \ \ \ \

 $\Box$

Computer simulations also confirms that the equilibrium point is not
stable.

\ \ \ \ \

In a general situation the above considerations exhibit the
existence of infinite oscillations for non-zero  coordinates.
Whether these oscillations may be cyclic? That is, whether the
cycles have a finite number of steps? Apparently  it is possible
under a certain starting connection between  values of all
coordinates. However,  the existence of exact finite cycles in the
"plus -- minus" model \ is the open question until now.

\ \ \ \ \ \

The  similar characteristic behavior has the conflict dynamical
system with the bilateral interaction of "minus -- plus" type -- the
model which describes the evolution of  biological species under
influence of some infection (the viral environment).

In such a case we have to define a vector of an initial
statistical distribution for a  virus infection along regions
$\Omega_i$:
 $${\bf q}=(q_1...,q_i,...,q_n), \ \
q_i:=\frac{Q_i}{Q}, $$ where
 $Q=
Q_1+\cdots +Q_n$. Then the evolution changes of $p^N_i$, $q^N_i$
 are governed  by the formulae
 similar to (\ref{pm}):
\begin{equation} \label{PQ} p_i^{N+1}=\frac{p_i^{N}(1-q
i^{N})}{z_{p,q}^{N}}, \ \ q_i^{N+1}=\frac{q_i^{N}(1+p
i^{N})}{z_{q,p}^{N}},\end{equation} where normalizing denominators
$z_{p,q}^{N}=1-({\bf p}^{N},{\bf q}^{N}), \ \ \ z_{q,p}^{N}=1+({\bf
q}^{N},{\bf p}^{N})$ ensure that  vectors ${\bf p}^N, {\bf q}^N$ are
stochastic. The opposite signs in the numerators of  (\ref{PQ}) have
now the following interpretation: minus means the decreasing of a
biological population  caused by infection, plus provides the growth
of  virus concentrations in regions with large local amount  of
biological species. In fact, according to  (\ref{PQ}),  the full
amount of  biological species and average virus concentrations are
stable because normative denominators  $z_{p,q}^{N}$, \
$z_{q,p}^{N}$ guarantee the non-annihilation of a population and the
natural dissipation for bacteria. Of course, the real quantitative
changes of these substances in various regions are determined by
(\ref{gl}) under the complex triple conflict interaction.

\ \ \ \ \ \ \

Let us now consider an abstract variant of the "minus -- minus" \
model. It may be interpreted as a situation of  conflict fighting
between couple of purely alternative opponents of type "infection
-- medicine". Shortly this alternative confrontation may be
written as "either -- or" \ that describes a tendency to exclusion
one to other from every region. In the terms of stochastic vectors
${\bf r}, {\bf q }$ the alternative interaction we represent by
the following  formulae:
\begin{equation}\label{MM}
q_i^{N+1}=\frac{q_i^{N}(1-r_i^{N})}{z^{N}}, \ \
r_i^{N+1}=\frac{r_i^{N}(1-q_i^{N})}{z^{N}},\end{equation} where a
value of  the normalizing  denominator $z^{N}=1-\theta^{N}$ strongly
dependents of the conflict index $\theta^{N}=({\bf q}^N,{\bf r}^N)$.

\ \ \ \ \ \

\begin{theorem} \label{Th2} ("Minus -- minus" \ model) \ Given a couple of stochastic vectors
  ${\bf q},{\bf r} \in {\mathbb R}^n_+, n>1$  assume the conflict index
$\theta=({\bf q},{\bf r})$  \  satisfies the inequalities:
$$0<\theta< 1 .$$ Then each trajectory of the conflict dynamical
system
 $$\{{\bf
q}^N,{\bf r}^N \} \stackrel{\divideontimes }{\longrightarrow }
\{{\bf q}^{N+1},{\bf r}^{N+1} \},  \ {\bf q}^{0}={\bf q}, {\bf
r}^{0}={\bf r}, \ N=0,1...$$ generated by  (\ref{MM}) goes with
necessity to the equilibrium state: $${\bf q}^{\infty}= \lim_{N
\to \infty}{\bf q}^{N}, \ \ {\bf r}^{\infty}= \lim_{N \to
\infty}{\bf r}^{N}$$ which is a fixed point: \ ${\bf q}^{\infty} =
{\bf q}^{\infty} { \divideontimes} {\bf r}^{\infty}, \ \ {\bf
r}^{\infty} = {\bf r}^{\infty} {\divideontimes} {\bf q}^{\infty}.$
Moreover, the limiting vectors are orthogonal \ $ {\bf q}^{\infty}
\perp   {\bf r}^{\infty},$ \  if  \ ${\bf q} \not= {\bf r}$, and
identical, \ $ {\bf q}^{\infty} = {\bf r}^{\infty}$, \  if the
initial vectors are equal ${\bf q} = {\bf r}$. In   latter case
the coordinates of limiting vectors are uniformly distributed:
$q_i^{\infty}=r_i^{\infty}=1/m$, \ where $m\leq n$ denotes the
amount  of non-zero initial  coordinates.
\end{theorem}

For the proof see  \cite{KoTC,KoTC1} and
\cite{KoDo,BKoDo,AKS,ABodK,BKS}.

We remark  that in  \cite{AKS} the formulae (\ref{MM}) was used
for the construction of  complex system which generalize the
well-known predator-prey model and is in fact some variant of the
vector  analog of Lotka-Volterra equations. Here we mention the
paper \cite{BT} where the idea of clusters (districts, regions)
was also used in three-dimensional discrete-time Lotka-Volterra
models.

Recall, that here we construct the discrete time models. However,
just below  we exhibit a pair formulae with continuous time:
$$\dot p(x,t)=\frac{p(\theta + r)}{m_{p,r} +\theta}, \ \ \dot
r(x,t)=\frac{r(\theta - p)}{m_{p,r}-\theta},$$ where   $p= p(x,t),
r=r(x,t), \ x\in \Omega$ denote the distribution densities of the
corresponding substances and \ $m_{p,r}=P\cdot R, \ \ \ \theta(t)=
\int_{\Omega} p(x,t)r(x,t) dx. $ We plan to study these equations in
consequent publications.

\section{The conflict triad}

We write  the  conflict dynamical system  that describes  the
evolution
 of simultaneously interacting triple
 substances    ${\bf P, R, Q}$ as follows:
\begin{equation}\label{ctr}\{{\bf  P}^{N}, {\bf R}^{N}, {\bf
Q}^{N} \} \stackrel{\divideontimes }{\longrightarrow } \{{\bf
P}^{N+1}, {\bf R}^{N+1}, {Q}^{N+1} \}, \ \ N=0,1...
\end{equation}
where the conflict map (composition) \ $\divideontimes$ is defined
below by formulae (\ref{gl}), (\ref{reg}). As above we assume that
substances ${\bf P, R, Q} $ have a common space of existence
$\Omega$, which is decomposed  in a natural way into the finite set
of separate regions, $\Omega=\bigcup_{i=1}^n\Omega_i, \ n\geq 2$.
The conflict triad is a complex system. It means that each substance
has an inner structure: $$ {\bf P}=(P_1...,P_n), \ {\bf
R}=(R_1...,R_n), \ {\bf Q}=(Q_1...,Q_n) $$ where the elements $P_i,
R_i, Q_i \ i=1. ,...,n$ determine the proper quantitative
description of the corresponding substances.

The  investigated substances of conflict triad have different
physical nature. That is why the concrete formulae of interactions
of every substance with  a complementary pair, namely  $ {\bf P } $
with the pair \ $\{{\bf  R, Q}\}$,  \ \ ${\bf R}$ with \ $\{{\bf P,
Q}\},$ and \ ${\bf Q}$ with \ \  $\{ {\bf P, R}\}$ \ are essentially
different one from another. We represent the complete mechanism of
interconnection that is contained in the conflict composition
$\divideontimes$ into two parts: \ formulae \ (\ref{gl}) which gives
the algorithm of quantitative changes of absolute values  $P_i, R_i,
Q_i$ in regions $\Omega_i$ and (\ref{reg}) which describes the
statistical law of redistribution of occupation  probabilities of
regions \ $\Omega_i$ \ by substances ${\bf P, R, Q} $.

The evolution of quantitative regional changes for   $P_i, R_i, Q_i$
is assigned  by the equations: $$P^{N+1}_i=\frac{P^{N}_i+d_1( R^{N}
_i-Q^{N}_i)}{Z^N_P}, \ \  $$ \begin{equation}\label{gl}
R^{N+1}_i=\frac{R^{N}_i+1/d_3\cdot Q^{N}_i/P^{N}_i}{Z^N_R},
\end{equation}
 $$Q^{N+1}_i=\frac{Q^{N}_i+d_2(R^{N}_i-Q^{N}_i)}{Z^N_Q},
 \ P_i^0=P_i \ R_i^0=R_i \ Q_i^0=Q_i \
N=0,1...,$$
 where \
parameters  $ \ d_1,d_2, d_3>0$  characterize the rate  of changes
intrinsic to the real model. The normalizing  denominators $Z^N_P,
Z^N_Q, Z^N_R$ ensure the stable global amount of the proper
substance in the whole space $\Omega$. Of course, the global amount
of every substance may be  additionally changed by virtue of
external circumstances, but  we do not consider such influence.
Here, for the sake of simplicity we assume that the  total
quantitative characteristics for each of substances are unchanged
and therefore we may write: $$P=\sum_i^n P_i=P^N, \ Q=\sum_i^n Q_i
=Q^N, \ R=\sum_i^nR_i=R^N,  \ N=0,1...$$

 We interpret formulae  (\ref{gl}) as
follows. The  quantitative growth of the biological species
 in   $\Omega_i$ region   on  $N+1$
step of the conflict fight is proportional to amount
  $P_i^{N}$  at  the previous moment of time  and to the difference
  values  with some
coefficient of the vital resource $ R_i^N$ and  the factor of
elimination threats $Q_i^N$. At some moment of time $N$ it may
happen that above difference has negative value. Then $P_i^N$ will
decrease quickly enough. However, such period of development has to
be   short. Otherwise, the system will be destroyed and loses its
physical sense, for example, if some of coordinates $P_i^N$,
$Q_i^N$, $ R_i^N$ becomes negative. Similarly, we interpret the
dependence of $Q_i^N$ from the same difference, but with another
coefficient. In turn the quantitative changes of  the vital resource
$R^{N}_i$ are very sensitive to the relative density of threat for
existence of biological population $Q^{N}_i/P^{N}_i$. In real models
the coefficient $d_3$  is small. We note that formally, according to
numerators in  formulae (\ref{gl}) all coordinates are increasing.
Nevertheless due to the normalizing denominators the dissipation
process courses. This automatically provides the decreasing of all
values $P^N_i, R_i^N, Q_i^N$ at each step of the conflict fight. Of
course, the concrete character of  interdependencies, their physical
interpretation, and the role of parameters is determined by  the
model of research.

The second part of our mechanism of the conflict interaction
$\divideontimes$ has purely probabilistic  inter-regional
character.
 To write it in the  mathematical terms it is necessary to transfer the
vectors $ {\bf P },$  \ ${\bf  R,  Q}$  into  stochastic ones:
$${\bf p}=(p_1...,p_n), \  \ {\bf r}=(r_1...,r_n), \  \ \ {\bf q
}=(q_1...,q_n), $$ where the coordinates $$p_i:=P_i/P,  \ r_i=R_i/R
\ \ q_i:=Q_i/Q \ \ i=1...,n $$ have a sense of probabilities to find
the corresponding   substance ${\bf P}$, ${\bf R}$ or ${\bf Q}$ in
$i$-th region.  In other words they are the occupation probabilities
 of  $\Omega_i$ by  $ {\bf P },$ \ ${\bf  R,  Q}$. In the course of
interactions the redistribution of these probabilities takes place.
The law of these changes is determined by the following formulae.
They are some  statistical variants of Lotka-Volterra discrete time
equations \cite{HoSi1,HoSi2} (c.f. with \cite{KoTC1}):
$$p^{N+1}_i=\frac{p^{N}_i(1+a(r^{N}_i-q^{N}_i))}{z^N_p}, $$
 \begin{equation}\label{reg} r^{N+1}_i=\frac{r^{N}_i(1-c
p^{N}_i-b q^{N}_i)}{z^N_r},\end{equation}
$$q^{N+1}_i=\frac{q^{N}_i(1+c^{-1} p^{N}_i-b r^{N}_i)}{z^N_q}, \
p^{0}_i=p_i, r^{0}_i=r_i, q^{0}_i=q_i, \ N=0,1,...,$$
     where parameters $a,b,c>0 $ characterize the intensity of the conflict
redistribution.

So the first of these formulae shows that the statistical
redistribution for the population substance ${\bf P}$ is maximal
in the region with the highest probability to find the vital
resource and the lowest threat to existence. In turn the
     second formula implies  the decreasing of a probability to find the
      vital resource
${\bf R}$ in a
 region where the biological population  is large
(because the later  utilizes the vital resource) and  there is a
high statistical  infection  concentration ${\bf Q}$. Finally,
according to the third formula in (\ref{reg}) the
   probability  of the infection threat (for the  population existence
   in  $i$-th region)  increases together with growing of
population and decreases under action of the
   vital resource as an alternative substance.
The denominators in   (\ref{reg}) \ provide that all vectors ${\bf
p}^N=(p^N_1,...,p^N_n), \  \ {\bf r}^N=(r^N_1,...,r^N_n), \ \ \
{\bf q }^N=(q^N_1,...,q^N_n), \ \ N=1,2,...$ have unite norms.

       To  complete the definition of the conflict map
       $\divideontimes$,  it is
necessary to fulfill the re-normalizing of the  vectors ${\bf
p}^{N+1}, \ \ {\bf r}^{N+1}, \ \ {\bf q }^{N+1}$ after using the
formulae (\ref{gl})  and  (\ref{reg}). It means transferring to the
final quantitative values of substances ${\bf P}$,  ${\bf Q}$ and
${\bf R}$ at each $(N+1)$-th step:
 \begin{equation}\label{ren} P^{N+1}_i=p_i^{N+1}\cdot P^{}, \
Q^{N+1}_i=q_i^{N+1}\cdot Q^{}, \ R^{N+1}_i=r_i^{N+1}\cdot R^{}.
\end{equation}

 Thus,  the   conflict composition as a map
$\divideontimes$ in  (\ref{ctr})  that  generates the  dynamical
system of conflict triad is entirely determined by formulas
(\ref{gl}), (\ref{reg}), (\ref{ren}).

In the present  work we made in fact only the first attempt to
construct and analyze   the simplest computer models of the
conflict triad. However even this activity finds out the series of
interesting observations usually  inherent to complex systems. In
particular, we establish the existence of fixed points (which are
 attractors), the existence of the stable limiting equilibrium  states, the appearance
  of cyclic orbits, which are attractors too,
  the critical bifurcation points, the oscillating trajectories
  without an evident  law of behavior, most similar to quasi-chaotic.
   So, we hope, the subsequent
research will lead to series of more deep results and useful
applications. \ \ \ \ \ \
\

 \begin{figure}[h] \centering
\includegraphics[height=10cm]{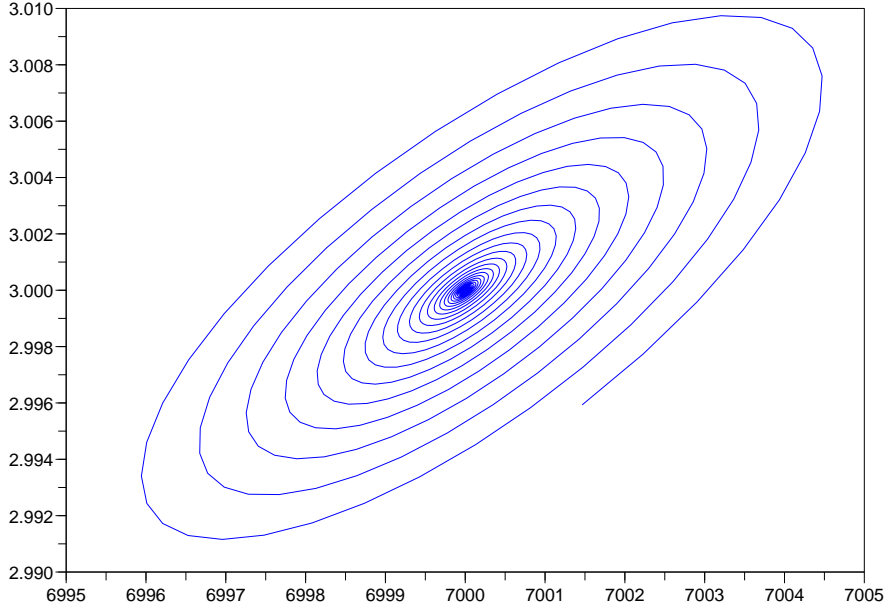}
\caption{{\it The equilibrium fixed point in the phase space}
$(P^N_1,Q^N_1)$}
\end{figure}

 \begin{figure}[h] \centering
\includegraphics[height=10cm]{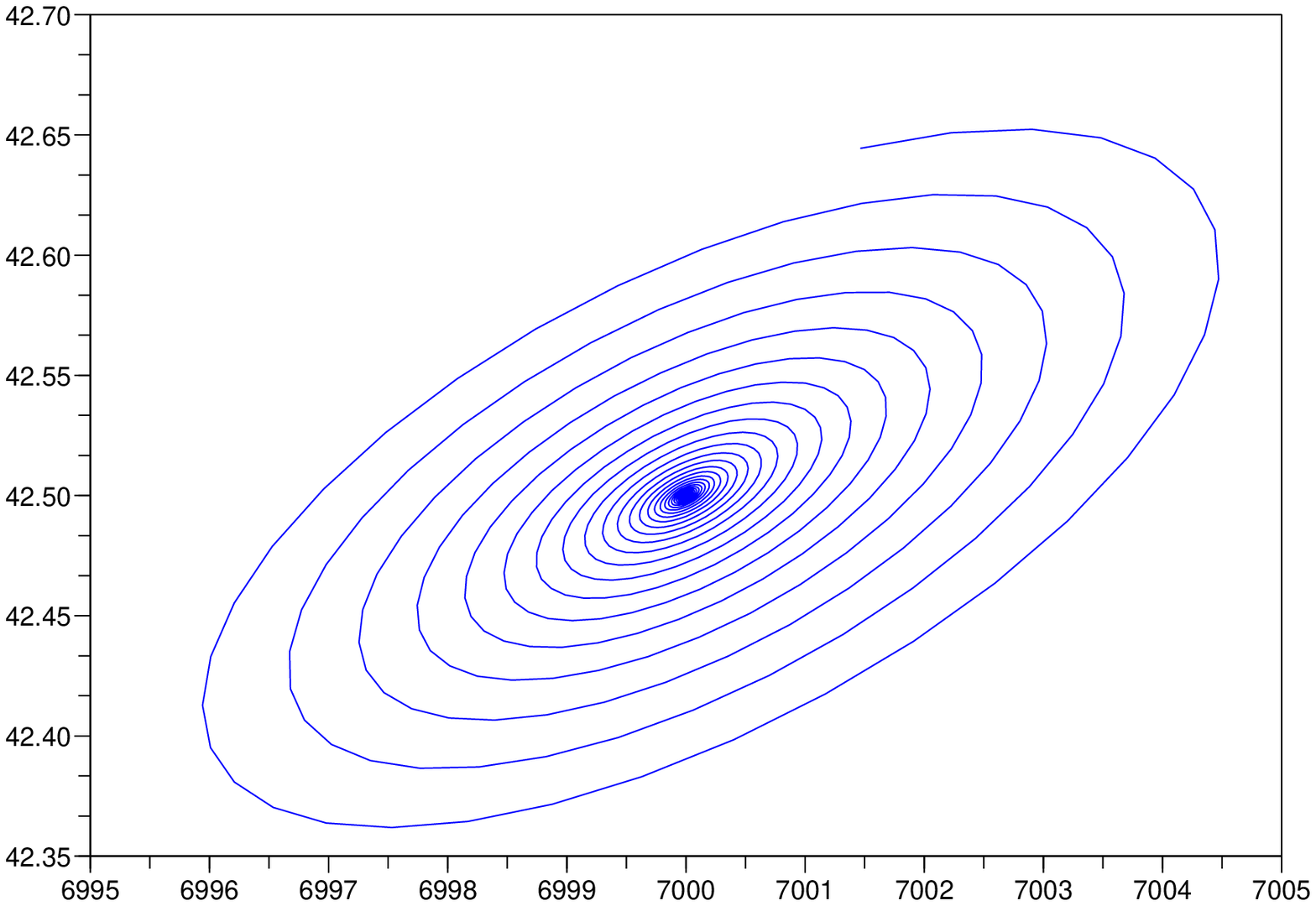}
\caption{{\it The equilibrium fixed point in the phase space}
$(P^N_1,R^N_1)$}
\end{figure}

\subsection{Computer models. The equilibrium state}

Due to formulae (\ref{gl}), (\ref{reg}), (\ref{ren}) the  conflict
triad is the highly  complex system. That is why the question of
existence of its state in which the appropriate substances are
situated in the equilibrium balance does not have an obvious answer.
In  physical reality we usually  observe that complex systems
possess such states as a rule. We mean the existence of the
dynamical equilibrium in a  complex  picture of the interaction
process between opponents with  contradict tendencies when the
various alternative substances coexist. We remark that the presence
of such state does not follow straightly from Theorems 3.1 and 3.2
ensuring the conditions for existence of the equilibrium states.
However one can conjecture about the possible mutual compensation of
oscillations which inherent to the models with bilateral
"plus-minus" \ interaction. It has to lead to stabilization and to
the equilibrium state \ as soon as the limiting compromise
distributions exist. Some analogs of such state are found in Theorem
3.2 for "minus--minus" \ models.

 \ \ \ \ \ \ \

\begin{theorem} \label{fp} The conflict triad dynamical system
(\ref{ctr}) defined by  (\ref{gl}) (\ref{reg}) (\ref{ren}), under
the condition  that all starting coordinates are non-zero,
possesses the equilibrium state. This state is determined by the
fixed point with coordinates $P_i^{\rm eq}, R_i^{\rm eq}, Q_i^{\rm
eq}, i=1...,n$ which are equal to the  arithmetic mean  values of
starting amounts of substances ${\bf P,R,Q}$ in regions
$\Omega_i$.
\end{theorem}
{\it Proof.} By virtue the assumption that global amounts $P, \ R,
\ Q$\ are constant,  the coordinates of any vector  ${\bf P, \ R,
\ Q}$\  cannot increase or decrease   simultaneously. So, if we
assume that at least one of coordinate is changed, for example,
becomes bigger, then there exists another one which will decrease
with necessity. However it is impossible since by (\ref{gl}),
(\ref{reg}), (\ref{ren}) these formulae are symmetric with respect
to permutations of indices and therefore all coordinates have the
same rights for changes. Thus, the state with
\begin{equation}\label{eq} P_i^{\rm eq}=1/n\sum_{k=1}^n P_k, \ \
R_i^{\rm eq}=1/n\sum_{k=1}^n R_k,  \ \  Q_i^{\rm
eq}=1/n\sum_{k=1}^n Q_k, \ \  i=1,...,n \end{equation} \ \ \ is
fixed. \ \ \ \ $\Box$

\ \ \ \ \ \

Surely the above Theorem \ref{fp} has the computer illustration.
If  in the concrete model one put the starting coordinates $P_i,
R_i, Q_i$ of vectors ${\bf P, R, Q}$ equal to middle-arithmetic
value
 $P_i^{\rm eq},
R_i^{\rm eq}, Q_i^{\rm eq}$, then they do not change  for any
$N\geq 1$. Thus, the corresponding state of the system is
equilibrium.

\begin{figure}[h] \centering
\includegraphics[height=10cm]{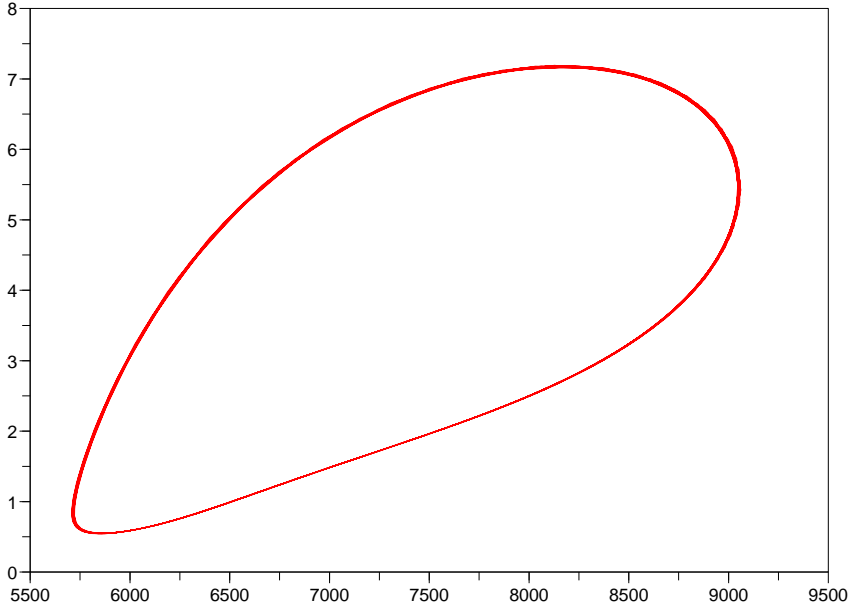}
\caption{{\it The cyclic attractor in the phase space}
$(P^N_1,Q^N_1)$}
\end{figure}

\begin{figure}[h] \centering
\includegraphics[height=10cm]{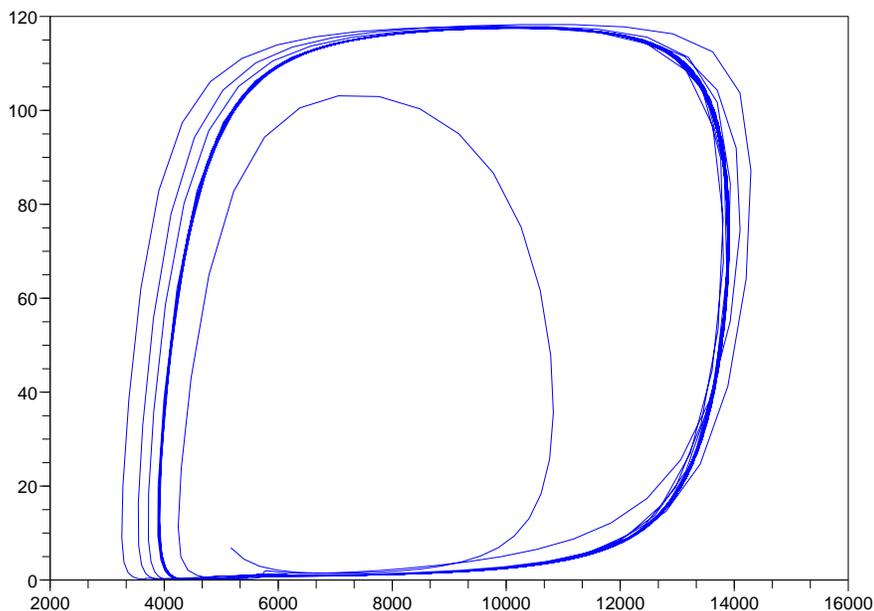}
\caption{{\it The almost square attractor in the phase space}
$(P^N_2,Q^N_2)$}
\end{figure}


\begin{figure}[h] \centering
\includegraphics[height=10cm]{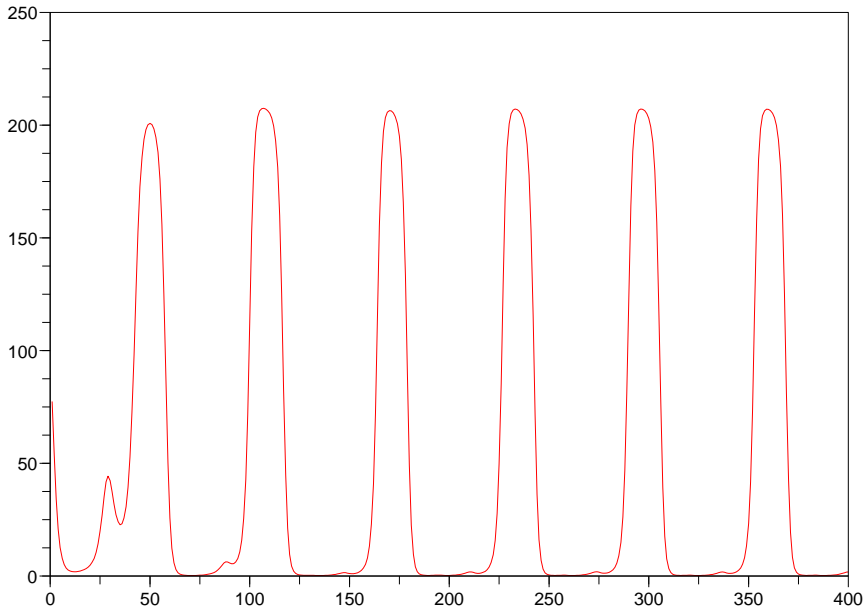}
\caption{{\it The periodic oscillation of the virus coordinate}
$Q_2^N$}
\end{figure}

\begin{figure}[h] \centering
\includegraphics[height=10cm]{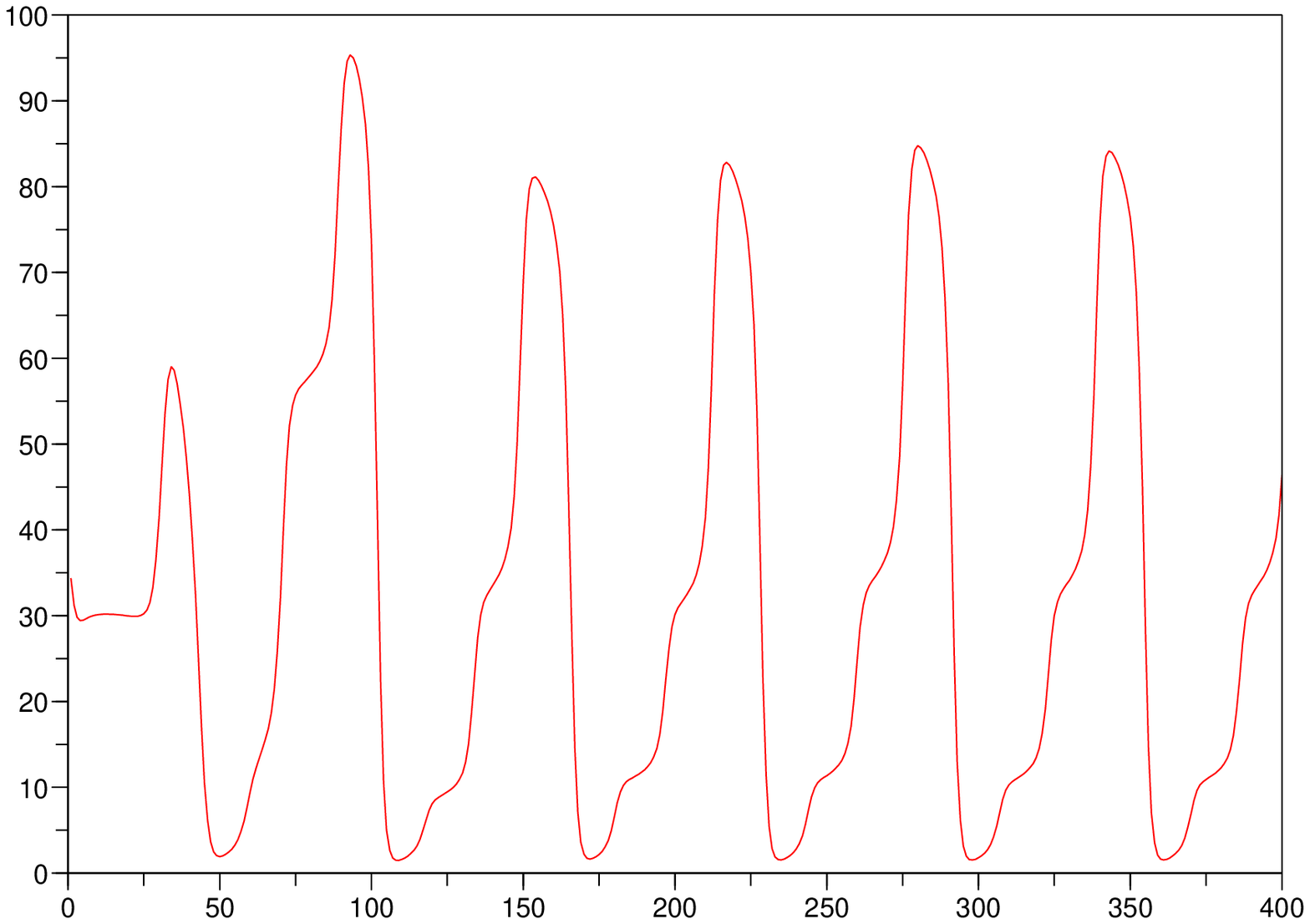}
\caption{{\it The periodic oscillation of the resource coordinate}
$Q_2^N$}
\end{figure}

\subsection{The equilibrium state is stable}

The fixed point from Theorem \ref{fp} is  in fact a one-point
attractor. This result we formulate as  follows

\begin{theorem} \label{opatr} The equilibrium state of the conflict triad
dynamical system  (\ref{ctr}) which is determined by the arithmetic
mean  values of the coordinates $P_i^{\rm eq}, R_i^{\rm eq},
Q_i^{\rm eq}, i=1...,n$ (see (\ref{eq})) is
 a one-point attractor.
\end{theorem}

\ \ \ \ \ \ \

Here we do not cite the  formal mathematical proof of this
theorem, but only the  computer illustration.

{\bf Example 1.} Consider the case with four regions, i.e., the
number of coordinates $n=4$. Let us choose the following meanings
 of parameters in formulas (\ref{gl}), (\ref{reg}): $$
d_1=d_3=0.09, \ d_2=0.01, \ a=0.1, \ b=0.6, \ c=0.1.$$ Note that
these values are guessed  "by hands". Under  their unreasonable
replacement, the model can "fall into pieces" \ (in particular, some
coordinates will go to infinite or become negative). Further, let
the population of biological species at initial moment of time has
the following quantitative distribution along regions:
$$P_1=9000; \ \ P_2=5000; \ \ P_3=2000; \ \ P_4=12000.$$ We
put the distributions of the vital resource and the infection
density to be essentially smaller: $$ R_1=30; \ \ R_2=80; \ \
R_3=50; \ \ R_4=10; \ \ Q_1=5; \ \
 Q_2=1; \ \ Q_3=2; \ \ Q_4=4. $$
The direct analysis of this computer model confirms Theorem
\ref{opatr}. We test the behavior of all coordinates to $N=2000 $
step. We observe, they demonstrate  the following features. Each
coordinate of all three vectors goes to an appropriate  value when
$N \to 2000 $. In the phase spaces we get the perfect spirals
which twist to the equilibrium point: $P_i^{\rm eq}, R_i^{\rm eq},
Q_i^{\rm eq}, i=1...,n$ with the arithmetic mean values  of the
initial dates.

It is important, that in the phase spaces $(P_i^N,R_i^N)$ and
$(P_i^N,Q_i^N)$ we get different spirals which twist in the opposite
directions to  points $(7000,42.5)$ \ and \ $(7000,3),$
respectively. We remark that  values  $P_i^N, R_i^N$ and $Q_i^N$
oscillate for a long time. They simultaneously approach  under $N
\to 2000$  to the fixed means $P_i^\infty=7000, \ R_i^\infty=42.5$
and $Q_i^N=3$ \ which exactly  are the middle values by regions of
the starting distributions of the proper substances (see Fig. 1 --
2).

\ \ \ \ \ \ \ \

Considering  of numerous examples, in particular, varying the
values of the initial coordinates demonstrate the stable character
 of the equilibrium state. No doubt,  the state of the conflict triad
  defined by the middle-arithmetic values is  attractive for
all trajectories close to this state. Thus, the proper point is a
local attractor. Moreover,  even large enough
 change of coordinates for the vector of biological populations
(for example, the replacement  $P_1=9000$ by $P_1=900$)
  does not destroy the attracting property  to  the equilibrium state.

\subsection{Cyclic attractors}

The equilibrium state from the previous example (the locally stable
fixed point) is not a global attractor. In particular, this state is
vanished  under the large enough change of single coordinate of the
vector that corresponds to an epidemic infection. Of course, under
any small variations of  coordinates,
 the system does not leave  the  attracting  phase to the
equilibrium state. Nevertheless, the replacement $ Q_2=1$ by $
Q_2=4$ transforms our system into another behavior phase. We
observed, in particular, the phase of attraction to the cyclic
orbit. Under the change  mentioned above all coordinates are not
attracted to a fixed point (the equilibrium state) but approximate a
cyclic trajectory. For more details, let us analyze the following
example.

 \ \ \ \ \

{\bf Example 2.} The loss  of the equilibrium state (a stable
fixed point) may take place under  some changing of the initial
coordinates and even only one of them.

 In particular, even changing  of the initial
values of any  substances. So, replacing $ R_2=80$ by $ R_2=40$ we
obtain the appearance of the cyclic attractor in the phase spaces
$(P_i^N,R_i^N)$ and $(P_i^N,Q_i^N)$ \ (see Fig. 3). That is, they
are achieved quickly enough, already at $N =400$ step. It is
interesting that  cyclic oscillations of
  $P_i^N, R_i^N$ and $Q_i^N$
take place around  points which are shifted with respect to the
arithmetic mean values of the initial coordinates (see the
previous example).

How  to explain this shift? There are also additional questions. For
example, why a decrease of full mass of vital resource  \ (from 170
to 130) leads also to passing of the system into the new phase? It
is not attracted  already by the fixed point. The evolution
trajectory of the conflict triad approximates the cyclic orbit. They
are attracted rather quickly  to a cycle of the egg-like form. That
is, in the limit all trajectories oscillate around fixed points
which are shifted in comparison with the initial mean values. All
these questions are open problems, although in \cite{AKS} we made
attempts to give some interpretations to the phenomena mentioned
above.

\ \ \ \ \ \ \

{\bf Example 3.} The exponential increasing  of the initial
coordinates $Q_i$  causes the appearance of almost square
attractor in the phase space. That is,
  $Q_i$ oscillates from zero to the maximal  value, nearly 120. In particular,
  the replacement
  $Q_2=10$ by $Q_2=100$ induced  the   oscillations of large amplitude, up to 200!
 \ (See Fig. 4).  In the same time the coordinates of resource vector
 $ R_i^N $ have periodic oscillations without achievement of the absolute maximum.

These observations  confirm an interesting practical effect. The \
substance which  corresponds to negative threats (an epidemic
infection) has large influence to the behavior of biological
population. So, the relatively small increase of initial values of
coordinates $Q_i$ substantially multiplies the negative effect on
other substances. Thus, the increasing of threats presses \ not
only to the existence of  biological population, but also at the
resource environment. In particular, it makes impossible the
achievement of  maximal values by them \ (See Fig. 5,6).

\subsection{The wave  of cyclic attractors}

 This phenomenon arises up under  a certain  increasing of
parameters $d_1, d_2, d_3$. We represent it at Fig. 7 where few
cycles are observed and they imposed one at another. These cycles
have different periods. In this case the trajectories  are
approximated consecutively to one of cyclic attractors, but only
for some time. We call this picture the wave of cyclic attractors.
Note that similar  behavior of complex biological dynamical
systems appear in practical situations, for example, in a case
with several epidemic sources.

The cyclic attractors for orbits can have not only egg-like shape,
but also considerably more difficult geometrical structures (see
Fig. 8). In particular, the cyclic attractor similar to the
${Carno}$-cycle (see Fig. 9) describes the behavior of a pair
"virus-resource".  It has close analogy with  work of a move
aggregate that consumes a certain resource, but survives resistance
and have to come back at the starting position.

\subsection{The quasi-chaotic behavior}

At the computer model that corresponds to  Fig. 8 we observe an
obvious evolution non-balanced for biological population and
infection, we call this dynamical phase as the quasi-chaotic
behavior. Thus, any regularity is  absent here and one can not find
even a slightly noticeable low in
 the behavior. We may find some
analogy with a chronic hidden  disease, that is healed, but not
cured. It increases non-periodically but does not reach a critical
stage in an organism which passed some threshold in its development.
In this situation it is impossible to go back to the state of
attracting to a stable point of equilibrium without external
influencing \ (See Fig. 9).

More deep analysis shows  the presence of bifurcation points in some
zones of phase space. Depending on parameter meanings $d_1,d_2,d_3$
and $a,b,c$ and  initial values of coordinates, the evolution of the
whole system passing this point may sharply change its direction. In
particular,  it is not structurally stable and can already not reach
the equilibrium state.

\begin{figure}[h] \centering
\includegraphics[height=10cm]{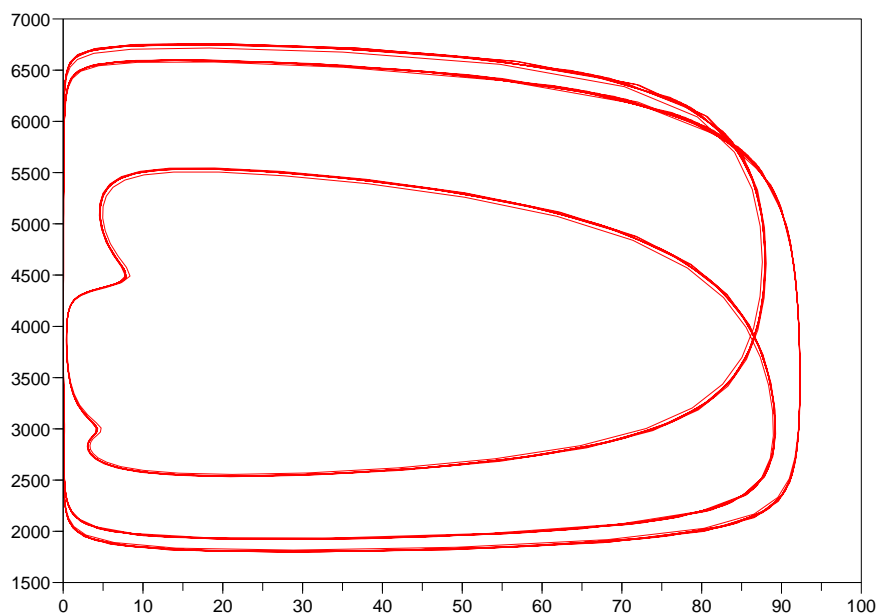}
\caption{{\it The multi-dimensional   attractor in the
phase-space} $(P^N_4,Q^N_4)$ \ \ ($d_1=0.001; d_2=0.000001;
d_3=0.0012; a=0.1; b=0.6; c=0.1; P_1=9000; P_2=5000; P_3=2000;
P_4=5; Q_1=50; Q_2=2; Q_3=1; Q_4=40; R_1=30; R_2=40; R_3=50;
R_4=10$)}
\end{figure}

\begin{figure}[h] \centering
\includegraphics[height=10cm]{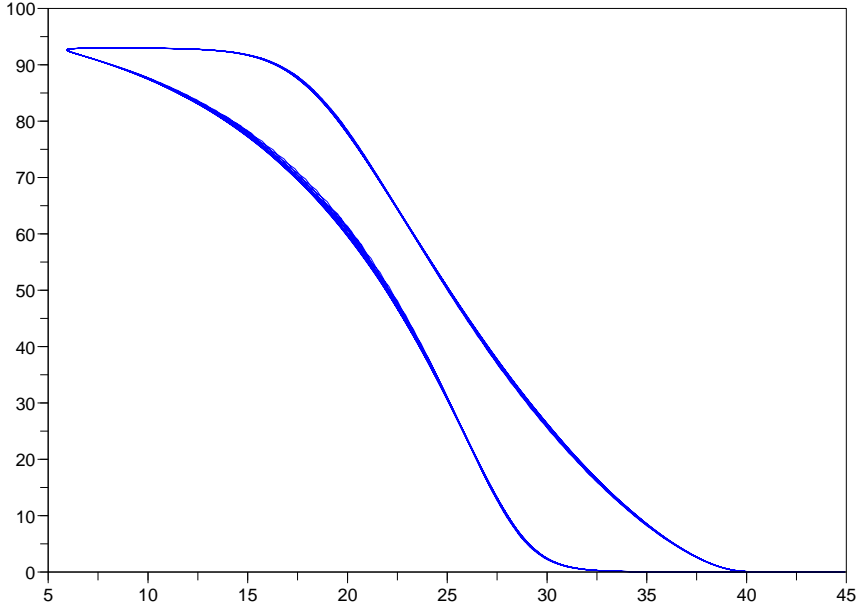}
\caption{{\it A Carno-type cycle in the phase-space}
$(R^N_2,Q^N_2)$ ({\it with} $d_3=0.0017$)}
\end{figure}

\begin{figure}[h] \centering
\includegraphics[height=10cm]{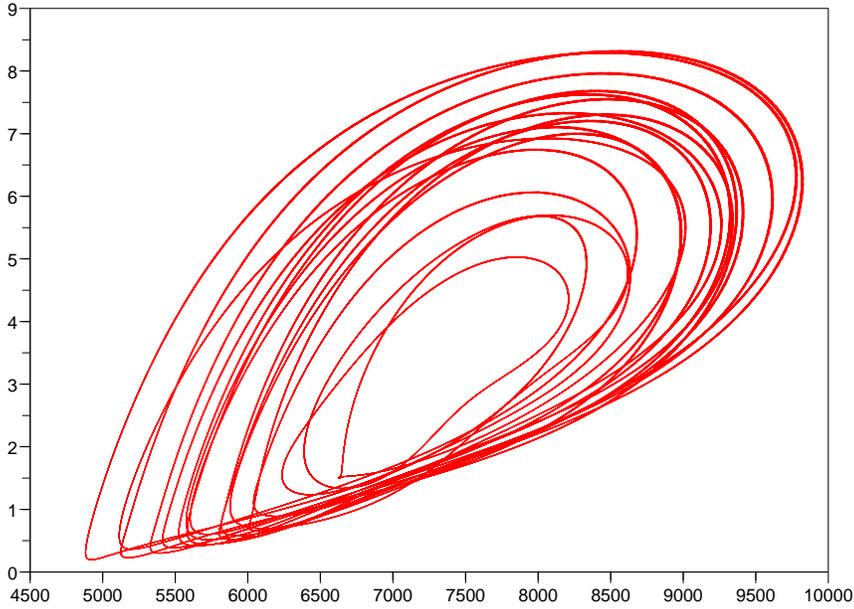}
\caption{{\it Several attracting sets with overlapping basins for}
$(P^N_1,Q^N_1)$ \ ($d_1=0.95; d_2=d_3=0.01; a=0.1; b=0.6; c=0.1;
P_1=9000; P_2=5000; P_3=2000; P_4=12000; Q_1=5; Q_2=1; Q_3=2;
Q_4=4; R_1=30; R_2=40; R_3=50; R_4=10$)}

\end{figure}

\section{Discussion}

The main result of our work consists in the construction of
complex models which demonstrate the presence of all basic phases
for coexistence of typical conflict triad substances: the state of
dynamical equilibrium (a stable fixed point), the phase of the
limit cyclic attractors (trajectories with periodically pulsating
evolution for each of substances), the area of bifurcation points,
that corresponds, in particular, to series  of the extremely
complex behaviors with transitions between some cascade of cyclic
attractors, and finally,  the phase close to chaotic which we call
the quasi-chaotic behavior. All these phases for coexistence of
conflict substances one can  meet rather  often in reality. So,
the state of dynamical equilibrium is typical  for the
simultaneous coexistence of different kinds  of bacteria inside a
living creature, when the highly dangerous microbes exist in a
healthy organism but do not exceed a critical concentration and
cause obvious rejections (that is an ordinary dynamical
equilibrium). Although the considerable enough external impact is
able to violate such equilibrium state and to result in  creation
of priorities to some  kind of bacteria and symptoms of a certain
illness. Other typical phase, the cyclic attractors, observed, for
example, when the season (cyclic) flu epidemics are caused by
exhausting of the valuable vital resources and by a too much
concentration (an overpopulation) of biological species in some
regions. This phase is character  of most oscillating processes in
our reality.

Certainly, the phases and states of the dynamical systems of
conflict triad mentioned  above appear  only under appropriated
values of parameters and starting coordinates. In the majority of
situations we observe an extremely complex, sometimes practically
chaotic behavior of our system. This all means that the theory
developed here is in fact only an introduction to study of the
conflict triad by the regional approach method. Of course, it can be
applied as a research  instrument for concrete tasks in conflicts
with triple opponents.

It is important that our model   in concrete settings allows to
determine the probability of infection for considered biological
species by an epidemic disease  in each separate region of common
space of existence. We recall that a key  point of our constructions
is a natural division  of the whole territory of existence $\Omega$
onto a set of regions $\Omega_i, \ i=1...,n$ with studying  of  the
local process of conflict interactions. Thus, the complete picture
of dynamical changes takes into account redistribution and migration
processes between regions, both for biological species and for
infections. In essence, a model gives not only the statistical
picture of distributions  but also the quantitative characteristics
for each of interacting substances in regions at discrete time. It
is possible to provide the prediction and information about the
power of infection risk, periods of relative safety and picks of
sharp growth of the infection density (epidemic) using the model.
Besides, it is also possible to pick regions with relative stability
or opposite ones with biggest growth of disease.

Finally we remark that it is possible to watch for sizes of
relations such as $R_i^N/P^N_i, Q^N_i/P^N_i$ which make sense of
distribution densities of positive vital resources and threats of
infection for existence of biological species  in $\Omega_i$ region.
These sizes are important  at the decision making  problem for safe
of existence of biological species  in a fixed region. The error of
uncertainty  of such decisions  depends on the range of changes for
parameters $d_1, d_2, d_3, a, b, c$ which are controlled by external
factors with respect to the system.

We hope that the dynamical model of conflict interaction between
triple intrinsic  elements offered here  may be used as a flexible
tool in  the problem of practical forecasts. Among important
dynamical parameters one can take such ones: the speed of
reproduction and spread of infection, their density and the local
concentration, the  rates of migration between regions and so on.

 \ \ \ \ \ \

\end{document}